\documentclass[review]{elsarticle}
\usepackage[utf8]{inputenc}
\usepackage{amsmath}
\usepackage{amsthm}
\usepackage{amssymb}
\usepackage{mathrsfs}
\usepackage{textcomp}
\usepackage{xcolor}
\usepackage{ulem}

\newtheorem{thm}{Theorem}[section]

\theoremstyle{definition}
\newtheorem{defn}{Definition}[section]

\theoremstyle{remark}
\newtheorem{rmk}{Remark}[section]

\journal{Journal of \LaTeX\ Templates}

\begin{document}

\begin{frontmatter}
\title{Global exponential stability and Input-to-State Stability of semilinear hyperbolic systems for the $L^{2}$ norm}
\author[a,b]{Amaury Hayat}
\address[a]{CERMICS, Ecole des Ponts ParisTech, 6-8 Avenue Blaise Pascal, Champs-sur-Marne, France.}
\ead{amaury.hayat@enpc.fr}
\address[b]{Department of Mathematical Sciences and Center for Computational and Integrative Biology, Rutgers University–Camden, 303 Cooper St, Camden, NJ, USA.}
\begin{abstract}
In this paper we study the global exponential stability in the $L^{2}$ norm of semilinear $1$-$d$ hyperbolic systems on a bounded domain, when the source term and the nonlinear boundary conditions are Lipschitz. We exhibit two sufficient stability conditions: an internal condition and a boundary condition.
This result holds also when the source term is nonlocal. Finally, we show its robustness by extending it to global Input-to State Stability in the $L^{2}$ norm with respect to both interior and boundary disturbances.
\end{abstract}
\begin{keyword}
Global stabilization \sep exponential stability \sep Lyapunov \sep hyperbolic systems; nonlinear \sep nonlocal \sep inhomogeneous
\MSC[2010] 35F60, 35F61, 93D09, 93D15, 93D20, 
93D30
\end{keyword}
\end{frontmatter}


\section{Introduction}

Hyperbolic systems can be found everywhere in sciences and nature: from biology \cite{Perthame2005}, to fluid mechanics, population dynamics \cite{TadmorTan}, electromagnetism, networks \cite{LeugeringSchmidt,DGLnet2010,EGMP} etc.
For this reason, they are of large importance for practical applications and the question of their stability and stabilization is paramount. For linear $1$-$d$ systems, studying the exponential stability or the stabilization can be achieved by looking at the eigenvalues and using spectral mapping theorems \cite{Lichtner,Renardy}. For nonlinear systems, the situation is much more tricky. 
For nonlinear systems, the situation is much more tricky. \textcolor{black}{In general the stabilities in different norms are not equivalent \cite{coron2015Nguyen}. \textcolor{black}{Indeed, for the same system, stabilities in different norms can require different criteria.}
For semilinear systems the spectral tools may still work (in contrast with quasilinear systems), but the resulting exponential stability may only hold locally, meaning for small enough perturbations.
 Worse, 
 most of the time spectral tools are hard to use when the system is inhomogeneous.}
Several tools were developed to deal with this situation and obtain local exponential stability
results. A first method is the characteristic analysis, which was originally used by Li and Greenberg in 1984 in \cite{Li1984} then generalized in \cite{1994-Li-book,Qin,Zhao,Wang} for quasilinear homogeneous hyperbolic systems in the $C^{1}$ norm. A second method is the use of basic Lyapunov functions\footnote{see \cite[Definition 1.4.3]{theseAH} for a proper definition and \cite{HS} for an overview of this method}. This method was, \textcolor{black}{for instance,} applied in \cite[Chapter 6]{BastinCoron1D} for general semilinear systems in the $H^{1}$ norm and quasilinear systems in the $H^{2}$ norm, but also in many particular cases \cite{BC2017,LeugeringH2,HS,SVgeneral,gugatH1}.
This will be our approach in this article. A third method is the backstepping method, a very powerful tool originally designed for finite-dimensional systems, modified for PDEs using a Volterra transform in \cite{Krsticsurvey}, \footnote{see \cite{xiang} for more details} and then used in \cite{CK2013,hu2015backstepping} for quasilinear hyperbolic systems in the $H^{2}$ norm. \textcolor{black}{Such backstepping approach was also used to derive controllability \cite{CN19a, CN19b} or finite-time stabilization \cite{CN17, CN20a, CN20b} in both parabolic and hyperbolic settings.} Other results using a more general transform were then introduced \cite{Zhang, MGC2018}. The main drawback of this method is that it involves controls that are usually using full-state measurements and cannot take the simple form of output feedback controllers (see \eqref{bound1}). Therefore these controls might be less convenient for practical implementation. Although sometimes observers can be designed to tackle this issue \cite{Meglio2013Observer}. \textcolor{black}{Other methods exist, as for instance the study of stability based on time delay systems introduced in \cite{coron2015Nguyen} where the authors give criteria for exponential stability in the $W^{2,p}$ norms for any $p\geq 1$ (see also \cite{CMS16}).}\\

So far, the nonlinear stability results for hyperbolic systems have been obtained in the $H^{1}$ norm for semilinear systems and for the $H^{2}$ norm for quasilinear systems. The $H^{1}$ and $H^{2}$ norm enabling to bound the nonlinear terms of the source term and of the transport term respectively, using the Sobolev embeddings $H^{p}([0,L];\mathbb{R})\subset C^{p-1}([0,L];\mathbb{R})$, for $p\geq1$. Other results have been shown for the $C^{0}$ and $C^{1}$ norm \cite{C1,C1_22}. For weaker norms, such as the $L^{2}$ norm, \textcolor{black}{one is usually unable to derive any exponential stability result when the system is nonlinear.
}
However, in this paper we show that having a Lipschitz source term, \textcolor{black}{with some condition on the size of the source,}
is enough to 
obtain
 the 
exponential stability
 in the $L^{2}$ norm for semilinear systems. Besides, in contrast with 
most of the
previous analyses cited above, this result holds for a nonlocal source term. Nonlocal source terms
are found in many important phenomena as population dynamics, material sciences, flocking, \textcolor{black}{traffic flow} \cite{TadmorTan, BeresPerthame,Keimernonlocal}, 
and open the door to many potential applications.
Moreover, while all the above previous approaches were dealing with local exponential 
stability,
we obtain here global exponential 
stability.
Concerning semilinear systems with Lipschitz source terms, one should
highlight
 the work of \cite{DFP} where the authors study the exponential 
 stability
 in 
$C^{0}$
 norm of a semilinear system with a diagonal and Lipschitz source term, and saturating boundary conditions. 
They give a potentially
large explicit bound on the basin of attraction, and they prove in addition the well-posedness in $L^{2}$.\\
 
 Finally, we show that these results can be extended to a wider notion:
 the Input-to-State Stability (ISS). The ISS measures the resilience of the stability of a system when adding disturbances in the boundary conditions or in the source term \cite{Sontag1989, KarafyllisKrstic}. \textcolor{black}{These disturbances could have many origins such as actuator errors, quantized measurments, uncertainties of model parameters, etc. The ISS} is therefore a more relevant notion from an application perspective, 
and 
is also \textcolor{black}{paramount for designing observers}. 
While exponential stability 
of nonlinear hyperbolic systems has been studied for several decades now, fewer results are known concerning this wider notion of ISS. Until recently,  the most up-to-date results were given in \cite[Part II]{KarafyllisKrstic}, for $L^{p}$ norms, $p\in\mathbb{N}^{*}\cup\{+\infty\}$ (see also \cite{prieuriss} for instance for nonautonomous systems), 
 and recently several works have been providing quite good conditions by extending exponential stabilization results obtained through Lyapunov approach to ISS results under the same conditions \cite{PrieurFerrante,YorgisBanda,BCH}. These results suffer however the same \textcolor{black}{limitations as the exponential stabilization results they are generalizing: local validity and strong norms.} One can also refers to \cite{Dashkovskiy,DashkovskiyMironchenko,Lhachemi2019,SaintVenantPI} for other  
 ISS results on hyperbolic systems in particular cases, and to \cite{KarafyllisKrstic} for a more detailed review on ISS results for PDEs in general. This paper is organized as follows: in Section \ref{s2} we state some definitions and our main result, which is proven in Section \ref{s3} using a Lyapunov approach. The well-posedness and the extension to ISS are dealt with in the Appendix.
 
  \section{Statement of the problem and main results}
 \label{s2}
A semilinear hyperbolic system can always be written in the following way \cite{LiYu}:
\begin{equation}
\partial_{t}\mathbf{u}+\Lambda(x)\partial_{x}\mathbf{u}+B(\mathbf{u},x)=0,
\label{sys1}
\end{equation}
where $\mathbf{u}(t,x)\in\mathbb{R}^{n}$, $\Lambda(x)$
 is a diagonal matrix with non vanishing eigenvalues, \textcolor{black}{$\Lambda :x\rightarrow \Lambda(x)$ belongs to $C^{1}([0,L])$} and $B\in C^{0}(L^{2}(0,L)\times [0,L],L^{2}(0,L))$ is the nonlinear source term, with $B(\mathbf{0},x)=0$. Note that
$B$ could be potentially nonlocal at it takes a function as argument, thus $B(\mathbf{u},x)$ refers here to $B(\mathbf{u}(t,\cdot),x)$. Throrough the article we will assume that $B(\cdot,x)$ is Lipschitz with respect to $\mathbf{u}$ with a Lipschitz constant $C_{B}$ in the following sense: for $\mathbf{u}$ and $\mathbf{v}$ two functions of $L^{2}(0,L)$,
\begin{equation}
\lVert B(\mathbf{u},\cdot)-B(\mathbf{v},\cdot)\rVert_{L^{2}}\leq C_{B}\lVert \mathbf{u}-\mathbf{v}\rVert_{L^{2}}.
\label{lipB}
\end{equation}
Of course, this assumption is satisfied if $B$ is local, takes argument in $\mathbb{R}^{n}\times[0,L]$ and  is Lipschitz with respect to the first argument, with a Lipschitz constant that might depend on $x$ but as a $L^{2}$ function.
We will come back to this special case later on in Remark \ref{rmklocal}.
When the system is equipped with a control static and exerted at the boundaries, the boundary conditions can be written in the following way:
\begin{equation}
\begin{pmatrix}
\mathbf{u}_{+}(t,0)\\
\mathbf{u}_{-}(t,L)
\end{pmatrix}=G\begin{pmatrix}
\mathbf{u}_{+}(t,L)\\
\mathbf{u}_{-}(t,0)
\end{pmatrix},
\label{bound1}
\end{equation}
where $G$ is a continuous and Lipschitz function such that $G(0)=0$. The notation $\mathbf{u}_{+}$ is used to refer to the components of $\mathbf{u}$ corresponding to positive propagation speeds $\Lambda_{i}>0$, whereas the notation $\mathbf{u}_{-}$ is used to refer to the components corresponding to negative propagation speeds. In the following, we assume without loss of generality that $\Lambda_{i}>0$ for $i\in\{1,...,m\}$ and $\Lambda_{i}<0$ for $i\in\{m+1,...,n\}$. Note that the boundary conditions \eqref{bound1} are nonlinear. 
As $G$ is Lipschitz, all of its components are Lipschitz, which implies that there exists a matrix $K$ such that
for any $i\in\{1,...,n\}$,
\begin{equation}
\left|G_{i}\begin{pmatrix}
\mathbf{u}_{+}(t,L)\\
\mathbf{u}_{-}(t,0)
\end{pmatrix}\right|\leq \sum\limits_{j=1}^{m} K_{ij}|u_{j}(t,L)|+\sum\limits_{j=m+1}^{n} K_{ij}|u_{j}(t,0)|.
\label{defK}
\end{equation} 
\begin{rmk}[Choice of $K$]
Of course the matrix $K=C_{G}I$, where $I$ is the identity matrix and $C_{G}$ the Lipschitz constant of $G$ would work. \textcolor{black}{However, there might be other matrices $K$ satisfying \eqref{defK} and some 
could lead to potentially less restrictive conditions in Theorem \ref{th1} than the matrix $C_{G}I$}
(see \eqref{condauxbords} below).
\end{rmk}
\textcolor{black}{System \eqref{sys1}, \eqref{bound1} with \eqref{lipB}, \eqref{defK}} is well posed in $L^{2}$ in the following sense: 
\begin{thm}[Well posedness]
For any $T>0$ and any $\mathbf{u}_{0}\in L^{2}(0,L)$ the Cauchy problem \eqref{sys1}--\eqref{bound1}, with initial condition $\mathbf{u}(0,\cdot)=\mathbf{u}_{0}$ has a unique solution $\mathbf{u}\in C^{0}([0,T],L^{2}(0,L))$. Moreover, 
\begin{equation}
\lVert \mathbf{u}(t,\cdot)\rVert_{L^{2}}\leq C(T)\lVert \mathbf{u}_{0}\rVert_{L^{2}},\text{ }\forall\text{ }t\in[0,T],
\label{estimate}
\end{equation}
where $C(T)$ is a constant depending only on $T$.
\label{th0}
\end{thm}
This theorem is shown in the Appendix. Most of 
the proof is a subcase of a remarkable result in \cite[Theorem A.1]{DFP}, where the authors study the framework of saturating boundary conditions. The only differences are some slight changes in the estimates to deal with a nonlocal functional and a density argument.
 These changes are indicated in \ref{reg}, 
 together with a proper definition of the notion of weak solution to System \eqref{sys1}, 
\eqref{bound1}.

\begin{rmk}
As it could be expected, the well posedness also holds for more regular solutions. In particular for any $\mathbf{u}_{0}\in H^{1}(0,L)$ satisfying the \textcolor{black}{compatibility} conditions \textcolor{black}{given by} \eqref{bound1}, the Cauchy problem \eqref{sys1}, \eqref{bound1} with initial condition $\mathbf{u}(0,\cdot)=\mathbf{u}_{0}$ has a unique solution $\mathbf{u}\in C^{0}([0,T],H^{1}(0,L)) \cap C^{1}([0,T],L^{2}(0,L))$. This is also shown in \ref{reg}.
\label{rmk1}
\end{rmk}

Before stating our main result, we recall the definition of exponential stability for the $L^{2}$ norm.
\begin{defn}[Exponential stability]
We say that System \eqref{sys1}--\eqref{bound1} is exponentially stable for the $L^{2}$ norm with decay rate $\gamma$ and gain $C$ if there exists constants $\delta>0$, $\gamma>0$, and $C>0$ such that for any $T>0$ and 
$\mathbf{u}_{0}\in L^{2}(0,L)$ such that $\lVert \mathbf{u}_{0}\rVert_{L^{2}}\leq \delta$, the Cauchy problem \eqref{sys1}--\eqref{bound1} with initial condition $\mathbf{u}(0,\cdot)=\mathbf{u}_{0}$ has a unique solution $\mathbf{u}\in C^{0}([0,T],L^{2}(0,L))$ and
\begin{equation}
\lVert \mathbf{u}(t,\cdot)\rVert_{L^{2}}\leq Ce^{-\gamma t}\lVert \mathbf{u}_{0}\rVert_{L^{2}}.
\end{equation}
Moreover, if
\begin{equation}
\delta=+\infty,
\end{equation}
then the system is said globally exponentially stable.
\end{defn}
We can now state our main result.
\begin{thm}
Let a system \textcolor{black}{be} of the form \eqref{sys1}, \eqref{bound1}, where $\Lambda\in C^{1}([0,L])$ and $B$ is 
Lipschitz with respect to $\mathbf{u}$.
\textcolor{black}{ If there exist $K\in M_{n}(\mathbb{R})$ satisfying \eqref{defK}, 
$J\in C^{1}([0,L]; M_{n}(\mathbb{R}))$ where 
$J(x)$ is a diagonal matrix
with positive coefficients, and 
$M\in C^{0}([0,L];M_{n}(\mathbb{R}))$, such that the following conditions are satisfied}
\begin{enumerate}
\item (Interior condition)
\textcolor{black}{
\begin{equation}
-(J^{2}\Lambda )'+J^{2}M+M^{\textcolor{black}{\top}}J^{2} 
\label{m1}
\end{equation}
is positive definite and there exists $D\in C^{1}([0,L]; M_{n}(\mathbb{R}))$ where $D(x)$ is a diagonal matrix
with positive coefficients, such that
\begin{equation}
C_{g}<\frac{\lambda_{m}}{2\max_{i,x}(D_{i})\max_{i,x}(D_{i}J^{2}_{i})},
\label{condint}
\end{equation}
where $C_{g}$ is the Lipschitz constant of $\textcolor{black}{g:=} B-M$ and $\lambda_{m}$ denotes the smallest eigenvalue of 
\begin{equation}
-D(\textcolor{black}{J^{2}\Lambda })'D+DJ^{2}MD+DM^{\textcolor{black}{\top}}J^{2}D,
\label{m2}
\end{equation}}
\item (Boundary condition) the matrix
\begin{equation}
\begin{split}
&\begin{pmatrix}J^{2}_{+}(L)\Lambda_{+}(L)&0\\
0& J^{2}_{-}(0)|\Lambda_{-}(0)|\end{pmatrix}\\
-&K^{\textcolor{black}{\top}}\begin{pmatrix}J^{2}_{+}(0)\Lambda_{+}(0)&0 \\ 0&J^{2}_{-}(L)|\Lambda_{-}(L)|\end{pmatrix}K
\label{condauxbords}
\end{split}
\end{equation}
is positive semidefinite,
\end{enumerate}
then the system is globally exponentially stable for the $L^{2}$ norm. Moreover the gain is $\lVert J^{-1} \rVert_{L^{\infty}}\lVert J \rVert_{L^{\infty}}$ \textcolor{black}{and an admissible decay rate is $\lambda_{m}(2\max_{i,x}(D_{i}J^{2}_{i}))^{-1}-C_{g}\max_{i,x}(D_{i})$}
\label{th1}
\end{thm}
We prove this theorem in Section \ref{s3}. 
Note that \eqref{condint} does not involve directly the Lipschitz constant of $B$ but the Lipschitz constant of $g \textcolor{black}{= B-M}$, which is $B$ minus a linear part that can be chosen. Of course, the Lipschitz constant of $B$ would be suitable by setting $M=0$, but other choices of $M$ could lead to less restrictive conditions. Let us note that the apparent complexity of the interior condition aims at giving a good explicit computable bound on $C_{g}$ for practical applications: indeed finding the values of $\lambda_{m}$ can be numerically solved. Besides, \textcolor{black}{choosing $D=Id$ or $K = C_{G}I$ would also give a sufficient condition that is simpler to write, but the sufficient condition would be more restrictive.}
\begin{rmk}[Linear case]
When $B$ is a local and linear operator we recover the result found in \cite[Proposition 5.1]{BastinCoron1D} \textcolor{black}{(see also \cite{DBC12} when $B$ is in addition marginally diagonally stable).  Indeed,} we can choose $M=B$, then $g=0$ and the interior condition is reduced to the existence of $J$, diagonal matrix with positive coefficients
 such that $-(\Lambda J^{2})'+J^{2}M+M^{\textcolor{black}{\top}}J^{2} $ is positive definite. 
\end{rmk}
\begin{rmk}[Local case]
In the special case where the system is local, i.e. $B$ is a function on $\mathbb{R}\times[0,L]$ and $B(\mathbf{u},x)=B(\mathbf{u}(t,x),x)$, the condition \eqref{condint} of the previous theorem can be slightly improved as follows: assume that $B$ is Lipschitz with respect to the first variable with a Lipschitz constant $C(x)\in L^{2}(0,L)$, then for any matrix $M$, $g=B-M$ is also Lipschitz with respect to the first variable and we can denote again its Lipschitz constant by $C_{g}(x)\in L^{2}(0,L)$. Then, the interior condition \eqref{condint} in Theorem \ref{th1} can be replaced by 
\begin{equation}
C_{g}<\frac{\lambda_{m}(x)}{\max_{i}(J^{2}_{i})(x)}\text{  or  }C_{g}<\mu_{m}(x)\frac{\max_{i}(J_{i})(x)}{\inf_{i}(J_{i})(x)},
\end{equation}
where $\lambda_{m}(x)$ and $\mu_{m}(x)$ are the smallest eigenvalues at a given $x$ of the matrix given by \eqref{m1} and \eqref{m2} respectively.
\label{rmklocal}
\end{rmk}
 \subsection{Input-to-State Stability}
In fact, this result can be extended to a more general notion: the Input-to-State Stability (ISS). This notion is more relevant when looking at practical implications as it takes into account the external disturbances that can arise. When such disturbance arise, System \eqref{sys1}, \eqref{bound1} is replaced by
\begin{equation}
\begin{split}
\partial_{t}\mathbf{u}+\Lambda(x)\partial_{x}\mathbf{u}+B(\mathbf{u},x)+\mathbf{d}_{1}(t,x)=0,\\
\begin{pmatrix}
\mathbf{u}_{+}(t,0)\\
\mathbf{u}_{-}(t,L)
\end{pmatrix}=G\begin{pmatrix}
\mathbf{u}_{+}(t,L)\\
\mathbf{u}_{-}(t,0)
\end{pmatrix}+\mathbf{d}_{2}(t),
\end{split}
\label{sys1ISS}
\end{equation}
where $\mathbf{d}_{1}$ and $\mathbf{d}_{2}$ are respectively the distributed and boundary disturbances. We define the 
ISS
as follows:
\begin{defn}[Input-to-State Stability]
We say that System \eqref{sys1ISS} is strongly Input-to-State stable (or ISS) with fading memory for the $L^{2}$ norm if there exists positive constants $\delta>0$, $C_{1}>0$, $C_{2}>0$, $\gamma>0$, such that for any $T>0$ and any $\mathbf{u}_{0}\in L^{2}(0,L)$ with $\lVert \mathbf{u}_{0} \rVert_{L^{2}}\leq \delta$ and $\lVert \mathbf{d}_{1}\rVert_{L^{2}}+\lVert \mathbf{d}_{2}\rVert_{L^{2}}\leq\delta$, there exists a unique solution $\mathbf{u}\in C^{0}([0,T],L^{2}([0,L]))$ to System \eqref{sys1}, \eqref{bound1}, and
\begin{equation}
\begin{split}
\lVert \mathbf{u}(t,\cdot)\rVert_{L^{2}}&\leq C_{1}e^{-\gamma t}\lVert \mathbf{u}_{0}\rVert_{L^{2}}\\
&+C_{2}\left(\lVert e^{-\gamma (t-s)}\mathbf{d}_{1}(s,x)\rVert_{L^{2}((0,t)\times(0,L))}\right.\\
&+\left. \lVert e^{-\gamma (t-s)}\mathbf{d}_{2}(s)\rVert_{L^{2}(0,t)}\right), \text{  for any  }t\in[0,T].
\end{split}
\label{ISS}
\end{equation}
Moreover, if $\delta=+\infty$, then the system is said to be globally strongly 
\textcolor{black}{ISS}
 with fading memory.
\label{defISS}
\end{defn}
This defines a strong notion of ISS with \textcolor{black}{an exponentially} fading memory. 
\textcolor{black}{The fading memory comes from the $e^{-\gamma(t-s)}$ in the $L^{2}$ norms of $\mathbf{d}_{1}$ and $\mathbf{d}_{2}$. It means that the influence of the disturbances at a given time $s$ decreases exponentially with time. One could have chosen other and less restrictive fading factors (see \cite[Chapter 7]{KarafyllisKrstic} for a more complete description of ISS estimates with fading memory).
}
The constants $C_{1}$ and $C_{2}$ are called the gains of the ISS estimate. When such notion of ISS cannot be achieved, weaker notions exist and can be found for instance in \cite{MironchenkoWirth}. 
We have the following result, analogous to Theorem \ref{th1}
\begin{thm}
Let a system \textcolor{black}{be} of the form \eqref{sys1ISS} where $\Lambda\in C^{1}([0,L])$, $\mathbf{d}_{1}\in \textcolor{black}{L^{2}}((0,T)\times(0,L))$, $\mathbf{d}_{2}\in \textcolor{black}{H^{1}}([0,T])$ and $B$ 
\textcolor{black}{is}
Lipschitz with respect to $\mathbf{u}$. 
If the condition \eqref{condint} is satisfied and the matrix defined by \eqref{condauxbords} is positive definite,
then the system is globally strongly ISS with fading memory for the $L^{2}$ norm.
\label{thmISS}
\end{thm}
The proof of this theorem is very similar to the proof of Theorem \ref{th1}. The only difference being that the assumption on \eqref{condauxbords} has to be slightly stronger than in Theorem \ref{th1} (positive definite instead of positive semidefinite). A way to adapt the proof of Theorem \ref{th1} is given in \ref{adapt}. Besides, the gains can again be computed explicitly as a function of $K$, $B$ and $\Lambda$ (see \eqref{estimateISSgains}).

 \section{Exponential stability in the $L^{2}$ norm}
 \label{s3}
 \begin{proof}[Proof of Theorem \ref{th1}]
 Let a semilinear system be of the form \eqref{sys1}, \eqref{bound1} with $\Lambda\in C^{1}([0,L],M_{n}(\mathbb{R}))$ and $B$ being $L^{2}$ 
 with respect to $\mathbf{u}$
 with Lipschitz constant $C_{B}$. We will first show Theorem \ref{th1} for $H^{1}$ solutions and then recover it for $L^{2}$ solutions using a density argument. Let $T>0$, and let $\mathbf{u}_{0} \in H^{1}(0,L)$. From Theorem \ref{th0} and Remark \ref{rmk1}, there exists a unique solution $\mathbf{u}\in C^{0}([0,T], H^{1}(0,L))\cap C^{1}([0,T],L^{2}(0,L))$ associated to this initial condition. 
 Let us now define the following Lyapunov function candidate:
 \begin{equation}
 V(\mathbf{u})=\int_{0}^{L}(J(x)\mathbf{u}(t,x))^{\textcolor{black}{\top}}J(x)\mathbf{u}(t,x) dx,
 \label{defV}
 \end{equation}
 where $J=\text{diag}(J_{1},...,J_{n})\in C^{1}([0,L],\mathcal{D}_{n}^{+}(\mathbb{R}^{n}))$, where $\mathcal{D}_{n}^{+}$ is the space of diagonal matrices with positive coefficients. The function $V$ is well defined on $L^{2}(0,L)$ and equivalent to 
\textcolor{black}{$\lVert \mathbf{u}(t,\cdot)\rVert^{2}_{L^{2}}$, as}
 \begin{equation}
\lVert\mathbf{u}(t,\cdot)\rVert_{L^{2}}^{2}\lVert J^{-1}\rVert_{L^{\infty}}^{-2}\leq V(\mathbf{u})\leq \lVert J\rVert_{L^{\infty}}^{2}\lVert\mathbf{u}(t,\cdot)\rVert_{L^{2}}^{2}.
\label{equivV}
 \end{equation}
We would like to show that $V$ decreases exponentially quickly along $\mathbf{u}$.
\textcolor{black}{Before going any further, let us comment on the choice of the form of this Lyapunov function candidate. Functions of this type are sometimes called \textit{basic quadratic Lyapunov function} or \textit{basic Lyapunov function for the $L^{2}$ norm} because they can be seen as the simplest functional equivalent of the $L^{2}$ norm. A commonly used Lyapunov function candidate for hyperbolic systems of conservation laws has the form \eqref{defV} with $J(x) = \text{diag}(q_{i}e^{-\mu s_{i} x})$ where $s_{i}=1$ if $\Lambda_{i}>0$ and $s_{i} = -1$ if $\Lambda_{i}<0$ and $q_{i}$ and $\mu$ are positive constants to be chosen. In our case however, such function might not work. This is due to the inhomogeneity and this a phenomena that can be seen in balance laws in general \cite{bastin2011coron}. For instance, in \cite{HS} is found a basic quadratic Lyapunov function that exists for any length $L>0$ provided good boundary conditions, while this could not happen with a basic quadratic Lyapunov function made of exponential weights.}
As $\mathbf{u}\in C^{1}([0,T],L^{2}(0,L))$, $V(\mathbf{u}(t,\cdot))$ can be differentiated with time, \textcolor{black}{and} we have
\begin{equation}
\begin{split}
\frac{d V(\mathbf{u}(t,\cdot))}{dt}&=\int_{0}^{L}2\mathbf{u}^{\textcolor{black}{\top}}J^{2}\partial_{t}\mathbf{u} dx\\
&=-\int_{0}^{L}2\mathbf{u}^{\textcolor{black}{\top}}J^{2}\Lambda\partial_{x}\mathbf{u} dx-2\int_{0}^{L}\mathbf{u}^{\textcolor{black}{\top}}J^{2}B(\mathbf{u},x) dx\\
&=-\left[\mathbf{u}^{\textcolor{black}{\top}}J^{2}\Lambda\mathbf{u}\right]_{0}^{L}+\int_{0}^{L}\mathbf{u}^{\textcolor{black}{\top}}(J^{2}\Lambda)'\mathbf{u} dx\\
&-2\int_{0}^{L}\mathbf{u}^{\textcolor{black}{\top}}J^{2}B(\mathbf{u},x) dx.
\end{split}
\end{equation}
We used here that $J$ and $\Lambda$ commute as they are both diagonal.
Now, 
\textcolor{black}{let $M\in C^{0}([0,L],M_{n}(\mathbb{R}))$} to be selected later on and set $g(\mathbf{u},x)=B(\mathbf{u},x)-M(x)\mathbf{u}(t,x)$ which is 
\textcolor{black}{again
Lipschitz in $\mathbf{u}$ in the sense of \eqref{lipB}}. \textcolor{black}{We} have
\begin{equation}
\begin{split}
\frac{d V(\mathbf{u}(t,\cdot))}{dt}&=-\left[\mathbf{u}^{\textcolor{black}{\top}}J^{2}\Lambda\mathbf{u}\right]_{0}^{L}+\int_{0}^{L}\mathbf{u}^{\textcolor{black}{\top}}(J^{2}\Lambda)'\mathbf{u} dx\\
&-2\int_{0}^{L}\mathbf{u}^{\textcolor{black}{\top}}J^{2}M\mathbf{u} dx-2\int_{0}^{L}\mathbf{u}^{\textcolor{black}{\top}}J^{2}g(\mathbf{u},x) dx\\
&=-\left[\mathbf{u}^{\textcolor{black}{\top}}J^{2}\Lambda\mathbf{u}\right]_{0}^{L}\\
&-\int_{0}^{L}\mathbf{u}^{\textcolor{black}{\top}}[-(J^{2}\Lambda)'+J^{2}M+M^{\textcolor{black}{\top}}J^{2}]\mathbf{u} dx\\
&-2\int_{0}^{L}\mathbf{u}^{\textcolor{black}{\top}}J^{2}g(\mathbf{u},x) dx\\
\end{split}
\label{dV}
\end{equation} 
where we used that $\mathbf{u}^{\textcolor{black}{\top}}J^{2}M\mathbf{u}=\mathbf{u}^{\textcolor{black}{\top}}M^{\textcolor{black}{\top}}J^{2}\mathbf{u}$, as it is a scalar. Now, we set
\begin{equation}
\begin{split}
I_{2}:&=\left[\mathbf{u}^{\textcolor{black}{\top}}J^{2}\Lambda\mathbf{u}\right]_{0}^{L},\\
I_{3}:&=\int_{0}^{L}\mathbf{u}^{\textcolor{black}{\top}}[-(J^{2}\Lambda)'+J^{2}M+M^{\textcolor{black}{\top}}J^{2}]\mathbf{u} dx\\
&+2\int_{0}^{L}\mathbf{u}^{\textcolor{black}{\top}}J^{2}g(\mathbf{u},x) dx
\label{defI3I2}
\end{split}
\end{equation}
We would like to show that under assumptions 
1. and 2. of Theorem \ref{th1}, $I_{2}$ is a nonnegative definite quadratic form with respect to the boundary conditions, and $I_{3}\geq \mu \lVert \mathbf{u}\rVert_{L^{2}}$ where $\mu$ is a positive constant. We will show that this is exactly the point of Assumptions 1. and 2..
Let us start with $I_{2}$. \textcolor{black}{From \eqref{bound1},}
\begin{equation}
\begin{split}
I_{2}&=\begin{pmatrix}\mathbf{u}_{+}(t,L)\\\mathbf{u}_{-}(t,L)\end{pmatrix}^{\textcolor{black}{\top}}J^{2}(L)\Lambda(L)\begin{pmatrix}\mathbf{u}_{+}(t,L)\\\mathbf{u}_{-}(t,L)\end{pmatrix}\\
&-\begin{pmatrix}\mathbf{u}_{+}(t,0)\\\mathbf{u}_{-}(t,0)\end{pmatrix}J^{2}(0)\Lambda(0)\begin{pmatrix}\mathbf{u}_{+}(t,0)\\\mathbf{u}_{-}(t,0)\end{pmatrix}\\
&=\sum\limits_{i=1}^{m}J_{i}^{2}(L)\Lambda_{i}(L)u_{i}^{2}(L)-\sum\limits_{i=m+1}^{n}J_{i}^{2}(0)\Lambda_{i}(0)u_{i}(0)^{2}\\
&+\sum\limits_{i=m+1}^{n}J_{i}^{2}(L)\Lambda_{i}(L)\left(G_{i}\begin{pmatrix}\mathbf{u}_{+}(t,L)\\\mathbf{u}_{-}(t,0)\end{pmatrix}\right)^{2}\\
&-\sum\limits_{i=1}^{m}J_{i}^{2}(0)\Lambda_{i}(0)\left(G_{i}\begin{pmatrix}\mathbf{u}_{+}(t,L)\\\mathbf{u}_{-}(t,0)\end{pmatrix}\right)^{2}.
\end{split}
\end{equation}
We set $x_{i}:=0$ if $i\in\{1,...,m\}$ and $x_{i}:=L$ if $i\in\{m+1,...,n\}$. Then using that $\Lambda_{i}>0$ for  $i\in\{1,...,m\}$ and $\Lambda_{i}<0$ otherwise, and using \eqref{defK},
\begin{equation}
\begin{split}
I_{2}&=\sum\limits_{i=1}^{n}J_{i}^{2}(L-x_{i})|\Lambda_{i}(L-x_{i})|u_{i}^{2}(L-x_{i})\\
&-\sum\limits_{i=1}^{n}J_{i}^{2}(x_{i})|\Lambda_{i}(x_{i})|\left(G_{i}\begin{pmatrix}\mathbf{u}_{+}(t,\textcolor{black}{L-x_{i}})\\\mathbf{u}_{-}(t,\textcolor{black}{L-x_{i}})\end{pmatrix}\right)^{2}\\
&\geq \sum\limits_{i=1}^{n}J_{i}^{2}(L-x_{i})|\Lambda_{i}(L-x_{i})|u_{i}^{2}(L-x_{i})\\
&-\sum\limits_{i=1}^{n}J_{i}^{2}(x_{i})|\Lambda_{i}(x_{i})|\left(\sum\limits_{j=1}^{n} K_{ij} |u_{j}(t,L-x_{\textcolor{black}{j}})|\right)^{2},\\
\end{split}
\label{I2withoutd}
\end{equation}
This can be rewritten as
\begin{equation}
\begin{split}
&I_{2}\geq \mathbf{Y}^{\textcolor{black}{\top}}
N
\mathbf{Y},
\end{split}
\label{estimI20}
\end{equation}
where $\mathbf{Y}$ is a vector with components $Y_{i}=|u_{i}(t,L-x_{i})|$ and $N$ is given by
\begin{equation}
\begin{split}
N=&\left(\begin{smallmatrix}J^{2}_{+}(L)|\Lambda_{+}(L)|&0\\
0& J^{2}_{-}(0)|\Lambda_{-}(0)|\end{smallmatrix}\right)\\
&-K^{\textcolor{black}{\top}}\left(\begin{smallmatrix}J^{2}_{+}(0)|\Lambda_{+}(0)|&0 \\0& J^{2}_{-}(L)|\Lambda_{-}(L)|\end{smallmatrix}\right)K.
\end{split}
\label{defN}
\end{equation}
 From \eqref{condauxbords} the matrix 
 $N$
 is positive semidefinite, thus
\begin{equation}
I_{2}\geq 0.
\label{estimI2}
\end{equation}
Let us now deal with $I_{3}$. 
Assume 
that the 
condition \eqref{condint} holds. Then
\textcolor{black}{there exists $D\in C^{1}([0,L], M_{n}(\mathbb{R}))$ such that $D(x)$ is diagonal with positive coefficients for any $x\in[0,L]$. Thus
\begin{equation}
 -D(J^{2}\Lambda)'D+DJ^{2}MD+DM^{\textcolor{black}{\top}}J^{2}D
\end{equation}
}
is a symmetric and definite positive matrix and we denote by $\lambda_{m}$ its smallest eigenvalue on $[0,L]$. We have from \eqref{defI3I2}, using Cauchy-Schwarz inequality and using the fact that $g$ is Lipschitz with $\mathbf{u}$ and the fact that $g(\mathbf{0},x)=B(\mathbf{0},x)=0$,
\textcolor{black}{
\begin{equation}
\begin{split}
I_{3}&\geq\int_{0}^{L}(D^{-1}\mathbf{u})^{\textcolor{black}{\top}}[-D(J^{2}\Lambda)'D+DJ^{2}MD+DM^{\textcolor{black}{\top}}J^{2}D](D^{-1}\mathbf{u}) dx\\
&-2\left(\int_{0}^{L}|D^{-1}\mathbf{u}|^{2}dx\right)^{1/2}\left(\int_{0}^{L}|DJ^{2}g(\mathbf{u},x)|^{2} dx\right)^{1/2}\\
&\geq\int_{0}^{L}(D^{-1}\mathbf{u})^{\textcolor{black}{\top}}[-D(J^{2}\Lambda)'D+DJ^{2}MD+DM^{\textcolor{black}{\top}}J^{2}D](D^{-1}\mathbf{u}) dx\\
&-2\max\limits_{i,x}(D_{i}J^{2}_{i}(x))\left(\int_{0}^{L}|D^{-1}\mathbf{u}|^{2}dx\right)^{1/2}\left(\int_{0}^{L}|g(\mathbf{u},x)|^{2} dx\right)^{1/2}\\
&\geq\int_{0}^{L}(D^{-1}\mathbf{u})^{\textcolor{black}{\top}}[-D(J^{2}\Lambda)'D+DJ^{2}MD+DM^{\textcolor{black}{\top}}J^{2}D](D^{-1}\mathbf{u}) dx\\
&-2\max\limits_{i,x}(D_{i}J^{2}_{i}(x))\left(\int_{0}^{L}|D^{-1}\mathbf{u}|^{2}dx\right)^{1/2}C_{g}\left(\int_{0}^{L}|\mathbf{u}|^{2} dx\right)^{1/2}\\
&\geq\int_{0}^{L}(D^{-1}\mathbf{u})^{\textcolor{black}{\top}}[-D(J^{2}\Lambda)'D+DJ^{2}MD+DM^{\textcolor{black}{\top}}J^{2}D](D^{-1}\mathbf{u}) dx\\
&-2C_{g}\max\limits_{i,x}(D_{i}J^{2}_{i}(x))\max\limits_{i,x}(D_{i}(x))\left(\int_{0}^{L}|D^{-1}\mathbf{u}|^{2}dx\right)^{1/2}\\
&\geq \lambda_{m}\lVert D^{-1}\mathbf{u}\rVert_{L^{2}}^{2}-2C_{g}\max\limits_{i,x}(D_{i}J^{2}_{i}(x))\max\limits_{i,x}(D_{i}(x))\lVert D^{-1}\mathbf{u}\rVert_{L^{2}}^{2}.
\end{split}
\end{equation}
}
Therefore if \textcolor{black}{$C_{g}<\lambda_{m}/(2\max_{i,x}(D_{i}(x))\max_{i,x}(D_{i}(x)J^{2}_{i}(x)))$} 
\textcolor{black}{then}
\begin{equation}
I_{3}\geq \mu \textcolor{black}{\lVert D^{-1}\mathbf{u}\rVert_{L^{2}}},
\label{estimI31}
\end{equation}
\textcolor{black}{with  $\mu =  \lambda_{m}-2C_{g}\max\limits_{i,x}(D_{i}J^{2}_{i}(x))\max\limits_{i,x}(D_{i}(x))>0$.}
Thus 
\textcolor{black}{ from \eqref{equivV}, the positive definiteness of $D$, hence $D^{-1}$, \eqref{dV}, \eqref{estimI2}, and \eqref{estimI31}}
\textcolor{black}{we can set $\gamma = \mu (\max(D_{i}J_{i}^{2}))^{-1} >0$} such that
such that for any $t\in [0,T]$.
\begin{equation}
\frac{d V(\mathbf{u}(t,\cdot)}{dt}\leq -\gamma V,
\end{equation}
and therefore 
\begin{equation}
V(\mathbf{u}(t,\cdot))\leq V(\mathbf{u}(s,\cdot))e^{-\gamma (t-s)},\text{ }\forall\text{ }0\leq s\leq t \leq T.
\end{equation}
From \eqref{equivV}, this implies that
\begin{equation}
\lVert \mathbf{u}(t,\cdot)\rVert_{L^{2}}\leq \lVert J^{-1}\rVert_{L^{\infty}}\lVert J\rVert_{L^{\infty}}e^{-\frac{\gamma}{2}(t-s)}\lVert\mathbf{u}_{0}\rVert_{L^{2}},
\label{expstab}
\end{equation}
which is exactly the estimate wanted \textcolor{black}{with decay rate $\gamma/2$}. So far this estimate is only true for $H^{1}$ solutions. However, it only involves the $L^{2}$ norm. Thus, as the system is well-posed in $C^{0}([0,T],L^{2}(0,L))$ and $\lVert\cdot\rVert_{L^{\infty}((0,T);L^{2}(0,L))}$ is lower semicontinuous, the estimate \eqref{expstab} also hold for $L^{2}$ solutions by density (more details on this argument can be found in the proof of \cite[Lemma 4.2]{Burgers}).

 \end{proof}
  \section{Adapting the proof in the ISS case}
 \label{adapt}
 In this section we show how to adapt the proof of Theorem \ref{th1} to get Theorem \ref{thmISS}. 
 \begin{proof}
 Let us consider System \eqref{sys1ISS} and let $T>0$. Let $\mathbf{u}_{0}\in H^{1}(0,L)$ and $\mathbf{u}\in C^{1}([0,T],H^{1}(0,L))$ the associated solution.
 Then, defining $V$ as in \eqref{defV}, and differentiating along $\mathbf{u}$,
 we obtain as previously
 \begin{equation}
\frac{d V(\mathbf{u}(t,\cdot))}{dt} = -I_{2}-I_{3}-2\int_{0}^{L}\mathbf{u}^{\textcolor{black}{\top}}J^{2} \mathbf{d}_{1}dx,
 \end{equation}
where $I_{2}$ and $I_{3}$ are given by \eqref{defI3I2}. Thus, using Young's inequality
 \begin{equation}
\frac{d V(\mathbf{u}(t,\cdot))}{dt} = -I_{2}-I_{3}+\varepsilon_{0} V+\frac{\lVert J^{2}\rVert_{L^{\infty}}}{\varepsilon_{0}}\lVert\mathbf{d}_{1}(t,\cdot)\rVert_{L^{2}}^{2},
 \label{estimateappenV}
 \end{equation}
 where $\varepsilon_{0}>0$ and can be chosen.
 As previously, from \eqref{condint}, $I_{3}\geq \mu V$ where $\mu>0$. \textcolor{black}{Therefore, choosing $\varepsilon_{0} = \mu/2$, we have
 \begin{equation}
 -I_{3}+\varepsilon_{0} V\leq -\frac{\mu}{2}V.
 \end{equation}
 }
 Concerning $I_{2}$, if we denote by $I_{2,0}$ the quantity in the absence of disturbances \textcolor{black}{(i.e. the quantity given by the first equality of \eqref{I2withoutd})} we get
 \begin{equation}
 \begin{split}
 I_{2}&=I_{2,0}-\sum\limits_{i=1}^{n}J_{i}^{2}(x_{i})|\Lambda_{i}(x_{i})|\left(d_{2,i}^{2}+2d_{2,i}G_{i}\begin{pmatrix}\mathbf{u}_{+}(t,L)\\\mathbf{u}_{-}(t,0)\end{pmatrix}\right)\\
&\geq \mathbf{Y}^{\textcolor{black}{\top}}N\mathbf{Y}-\sum\limits_{i=1}^{n}J_{i}^{2}(x_{i})|\Lambda_{i}(x_{i})|\left(1+\frac{1}{\varepsilon}\right)d_{2,i}^{2}\\
&-\varepsilon\mathbf{Y}K^{\textcolor{black}{\top}}\begin{pmatrix}J^{2}_{+}(0)|\Lambda_{+}(0)|&0 \\0& J^{2}_{-}(L)|\Lambda_{-}(L)|\end{pmatrix}K\mathbf{Y},
\end{split}
\end{equation}
 \textcolor{black}{where we used Young's inequality and} where $N$ is the matrix given in \eqref{defN}, $\mathbf{Y}$ is defined as in \eqref{estimI20}, and $\varepsilon>0$ is 
 to be chosen. Using the definition of $N$ and the fact that $N$ is positive definite (and not positive semidefinite in contrast with Theorem \ref{th1}), we get by continuity that there exists $\varepsilon>0$ such that
\textcolor{black}{
 \begin{equation}
N-\varepsilon K^{\textcolor{black}{\top}}\begin{pmatrix}J^{2}_{+}(0)|\Lambda_{+}(0)|&0 \\0& J^{2}_{-}(L)|\Lambda_{-}(L)|\end{pmatrix}K\;\;\;\text{ is semipositive definite}.
 \end{equation}
 Therefore,
  $I_{2}\geq -(1+\varepsilon^{-1})\lVert J\rVert_{\infty}^{2}\lVert \Lambda\rVert_{\infty} |\mathbf{d}_{2}(t)|^{2}$ and
\eqref{estimateappenV} becomes}
 \begin{equation}
 \begin{split}
 \frac{dV(\mathbf{u}(t,\cdot))}{dt}\leq&-\frac{\mu}{2} V+\frac{2\lVert J\rVert_{L^{\infty}}^{2}}{\mu}\lVert\mathbf{d}_{1}(t,\cdot)\rVert_{L^{2}}^{2}\\
 &+(1+\varepsilon^{-1})\lVert J\rVert_{\infty}^{2}\lVert \Lambda\rVert_{\infty}  |\mathbf{d}_{2}(\textcolor{black}{s})|^{2},
 \end{split}
 \end{equation}
 thus, using Gronwall's Lemma, 
 \begin{equation}
 \begin{split}
&V(\mathbf{u}(t,\cdot))\leq V(\mathbf{u}_{0})e^{-\frac{\mu t}{2}}\\
&+\frac{2\lVert J\rVert_{L^{\infty}}^{2}}{\mu}\int_{0}^{t}e^{-\frac{\mu}{2} (t-s)}\left(\lVert\mathbf{d}_{1}(s,\cdot)\rVert_{L^{2}}^{2}\right.\\
&\left.+\frac{\mu}{2}(1+\varepsilon^{-1})\lVert \Lambda\rVert_{\infty} |\mathbf{d}_{2}(t)|^{2}\right)ds,
\end{split}
\label{estimateISSgains00}
 \end{equation}
which, together with \eqref{equivV} and the concavity of the square root function
 gives \textcolor{black}{
 \begin{equation}
  \begin{split}
&\lVert \mathbf{u}(t,\cdot)\rVert_{L^{2}}\leq \lVert J^{-1} \rVert_{L^{\infty}} \lVert J \rVert_{L^{\infty}}\lVert \mathbf{u}_{0}\rVert_{L^{2}} e^{-\frac{\mu t}{4}}\\
&+\lVert J^{-1} \rVert_{L^{\infty}}\lVert J\rVert_{L^{\infty}}\sqrt{\frac{2}{\mu}\max\left(1,\frac{\mu}{2}(1+\varepsilon^{-1})\lVert \Lambda\rVert_{L^{\infty}}\right)}\left(\lVert e^{-\frac{\mu}{2} (t-s)} \mathbf{d}_{1}(s,x)\rVert_{L^{2}((0,t)\times(0,L))}+\right.\\
&\left.+ \lVert \mathbf{d}_{2}(t) \rVert_{L^{2}(0,t)}\right),
\end{split}
\label{estimateISSgains}
 \end{equation}
which is the ISS estimate wanted and this holds for any $H^{1}$ solutions}. And, by density, this holds also for any $L^{2}$ solutions. Note that 
  the gains of the estimate can again be computed explicitly. This ends the proof of Theorem \ref{thmISS}
 \end{proof}
 \section{\textcolor{black}{Numerical simulations}}
 \label{examples}
\begin{figure}
{\centering \includegraphics[height = 0.3\textheight]{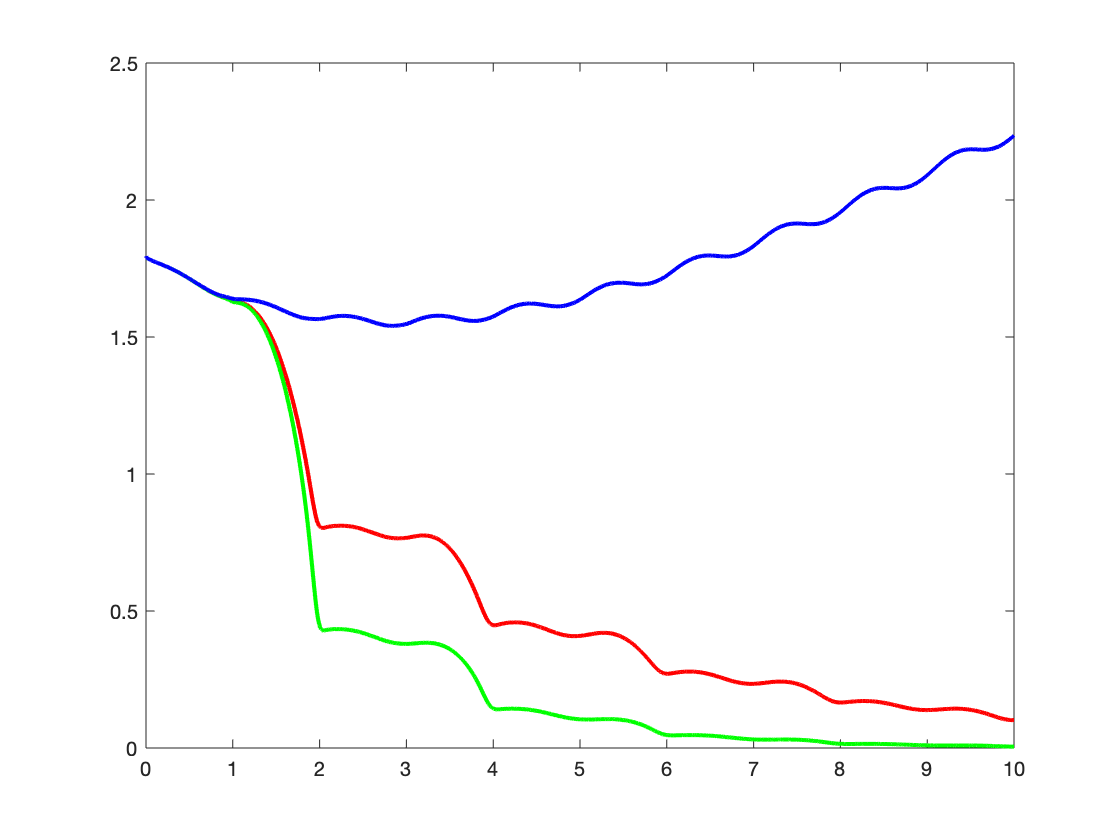} \par}
\caption{\textcolor{black}{Stability of the system \eqref{systemex} in open-loop (blue) and closed-loop with $k =3/4$ (green) and $k=1/2$ (red). horizontal axis represents time, and vertical axis represents the $L^{2}$ norm of the solution with initial condition $u_{1,0}(x) = \sqrt{2\pi x}$ and $u_{2,0}(x) = e^{-2\pi x}$. \label{fig1}}}
\end{figure}

\textcolor{black}{
In this section we present a numerical illustration of the previous result on a simple example. We consider a 
system inspired from \cite[Section 5.6]{BastinCoron1D} and given as
\begin{equation}
\begin{split}
\partial_{t}u_{1}+\partial_{x}u_{1} &= cL^{-1}\sin\left(\int_{0}^{L}u_{2}(t,x)dx\right)\\
\partial_{t}u_{2}-\partial_{x}u_{2} &= cL^{-1}\sin\left(\int_{0}^{L}u_{1}(t,x)dx\right)\\
u_{1}(t,0)-u_{2}(t,0) &= 0\\
u_{1}(t,L)-u_{2}(t,L) &= k u_{1}(t,L)\\
\end{split}
\label{systemex}
\end{equation}
where one boundary condition can be imposed through a design parameter $k$ while the other one is imposed.
Note first that in open-loop, i.e. $k=0$, the null steady-state is an unstable steady-states for any $c\in \mathbb{R}$ and any length of the domain $L>0$. Indeed, there is a continuum of travelling wave solutions: for any $\varepsilon>0$ 
\begin{equation}
\left\{
\begin{split}
u_{1}(x) &= \varepsilon e^{\frac{2\pi i}{L}(t-x)}\\
u_{2}(x) &= \varepsilon e^{-\frac{2\pi i}{L}(t-x)}
\end{split}
\right.
\end{equation}
is a solution of \eqref{systemex} with $k = 0$.
Nevertheless, Theorem \ref{th1} can be applied to find a feedback in closed loop as long as $|c| L <1/2$: set $M=0$, 
$D=Id$, $\varepsilon>0$ to be defined, and $k=\sqrt{1/(1+2L\varepsilon^{-1})}$.
Set also, $J = \text{diag}(\sqrt{L+\varepsilon- x}, \sqrt{L+\varepsilon+x})$, one has
$\Lambda =\text{diag}(1,-1)$ therefore  $-(J^{2}\Lambda)' = I_{d}$ and therefore is positive definite with smallest eigenvalue $1$. Besides $\max_{i,x}(J_{i}^{2}) = \varepsilon+2L$ and 
\begin{equation}
\begin{split}
\lVert g(U)-g(V)\rVert_{L^{2}}^{2}
=&\frac{1}{L^{2}}\int_{0}^{L}\left|\begin{pmatrix}c\sin\left(\int_{0}^{L}U_{2}(x)dx\right)\\
c\sin\left(\int_{0}^{L}U_{1}(x)dx\right)-
\end{pmatrix}-\begin{pmatrix}c\sin\left(\int_{0}^{L}V_{2}(x)dx\right)\\
c\sin\left(\int_{0}^{L}C_{1}(x)dx\right)-
\end{pmatrix}\right|^{2}dx\\
&\leq |c|L^{-1}\left[\left(\int_{0}^{L}|U_{2}-V_{2}|dx\right)^{2}+\left(\int_{0}^{L}|U_{1}-V_{1}|dx\right)^{2}\right]\\
&\leq |c|\lVert U-V\rVert_{L^{2}}^{2}.
\end{split}
\end{equation}
Hence, condition \eqref{condint} becomes $|c|<(\varepsilon+2L)^{-1}$. Now, as $|c|<(2L)^{-1}$, one can choose $\varepsilon=3(|c|^{-1}-2L)/4$ such that condition \eqref{condint} is satisfied. Finally, one can easily check that condition \eqref{condauxbords} becomes
\begin{equation}
(1-k)^{2}\leq \frac{\varepsilon}{\varepsilon+2L},
\end{equation}
which is also satisfied from our definition of $k$. Thus Theorem \ref{th1} applies and the system is globally stable for the $L^{2}$ norm. 
On Figure \ref{fig1} we represent the $L^{2}$ norm of the solution for various values of $k$ when $c = 1/4$ and $L=1$. In blue is represented the open-loop situation (i.e. $k=0$), in green the closed-loop situation with $k = 3/4$, and in red $k =1/2$.}

\section{\textcolor{black}{Conclusion and perspective}}
\textcolor{black}{We derived sufficient conditions for the global stability in the $L^{2}$ norm of semilinear systems with Lipschitz boundary conditions and source term (potentially nonlocal). We also showed that a strong ISS property with respect to boundary and internal disturbances holds globally under the same conditions. 
This result could have many applications in practice. Knowing whether such conditions are optimal for the existence of a basic quadratic Lyapunov function, at least for $n=2$ as it is in the linear and local case, is an open question. Another interesting direction for future works would be to try to extend, at least partially, these results to quasilinear but Lipschitz nonlocal systems.
}
 \section*{Acknowledgments}
The author would like to thank Jean-Michel Coron for his advice and many interesting discussions. The author would also like to thank Benedetto Piccoli for interesting discussions, the NSF for support via the CPS Synergy project "Smoothing Traffic via Energy-efficient Autonomous Driving" (STEAD) CNS 1837481, the French Corps des IPEF,  the ANR Finite4Sos, and the INRIA CAGE team.
 \appendix
 \section{Well-posedness of the system}
  \label{reg}
 In this section we deal with the well-posedness of the system and extend \cite[Theorem A.1]{DFP} to get Theorem \ref{th0}. But first, we 
give 
the definition of a weak $L^{2}$ solution for System \eqref{sys1},\eqref{bound1}.
 \begin{defn}
 \label{defsol}
Let $\mathbf{u}_{0}\in L^{2}(0,L)$. We say that $\mathbf{u}\in C^{0}([0,+\infty);L^{2}(0,L))$ is an $L^{2}$ solution of the Cauchy problem \eqref{sys1}, \eqref{bound1}, $\mathbf{u}(0,\cdot)=\mathbf{u}_{0}$, if for every $T>0$ 
there exists a sequence of functions $\mathbf{u}_{0,n}\in H^{1}(0,L)$ satisfying \eqref{bound1} and such that 
\begin{equation}
\begin{split}
&\mathbf{u}_{0,n}\rightarrow \mathbf{u}_{0}\text{ in }L^{2}(0,L),\\
&\mathbf{u}_{n}\rightarrow \mathbf{u}\text{ in }C^{0}([0,T],L^{2}(0,L)),
\end{split}
\label{convdef}
\end{equation}
where $\mathbf{u}_{n}\in C^{0}([0,T],H^{1}(0,L))$ is a weak solution of \eqref{sys1}, \eqref{bound1} with initial condition $\mathbf{u}_{0,n}$, i.e. $\mathbf{u}_{n}$ satisfies \eqref{bound1} and for 
any $\phi\in C^{1}([0,T];C^{1}_{c}((0,L);\mathbb{R}^{n}))$
we have
\begin{equation}
\begin{split}
\int_{0}^{L}&\int_{0}^{T}\partial_{t}\phi^{\textcolor{black}{\top}}\mathbf{u}_{n}+\partial_{x}\phi^{\textcolor{black}{\top}}\Lambda(x)\mathbf{u}_{n}\\
&+\phi^{\textcolor{black}{\top}}(\Lambda_{x}\mathbf{u}_{n}-B(\mathbf{u}_{n},x))dt\text{ }dx\\
=&
\int_{0}^{L}\left[\phi(\cdot,x)^{\textcolor{black}{\top}}\mathbf{u}_{n}(\cdot,x)\right]_{0}^{T} dx.
\end{split}
\label{L2sol}
\end{equation}
 \end{defn}
 \textcolor{black}{
 \begin{rmk}
 As noted in \cite{DFP}, this definition is slightly different from the definition given in \cite[Definition A.3]{BastinCoron1D} when looking at linear systems. The reason comes from the nonlinear boundary conditions which may prevent the adjoint of the boundary operator from existing.
 \textcolor{black}{Of course, in the linear case, a solution in the sense of \cite[Definition A.3]{BastinCoron1D} is also a solution in the sense of Definition \ref{defsol}.}
 \end{rmk}
 }
With this definition in mind, we prove Theorem \ref{th0}, by \textcolor{black}{slightly} adapting the proof of \cite[Theorem A.1]{DFP}
 \begin{proof}[Proof of Theorem \ref{th0}]
Let $T>0$. We define the operator $\mathcal{A}=-\Lambda(x)\partial_{x}$ on the domain $D(\mathcal{A})$ defined by
 \begin{equation}
 D(\mathcal{A})=\{\mathbf{u}\in H^{1}(0,L) | \mathbf{u} \text{ satisfies \eqref{bound1} }\} .
 \end{equation}
We also consider $B$ as an operator on the domain $D(B)=L^{2}(0,L)$, \textcolor{black}{and in the following $Bf$ refers to $B(f(\cdot),x)\in L^{2}(0,L)$}. Observe that $D(\mathcal{A}+B)=D(\mathcal{A})$. \textcolor{black}{First of all, we can restrict ourselves to the case where $\Lambda$ has only positive components. Indeed, if not, we define $\mathbf{v}=(v_{i})_{i\in\{1,...,n\}}$ by
\begin{equation}
\left\{\begin{split}
v_{i}(t,\cdot)&=u_{i}(t,\cdot)\text{ if }i\in\{1,...,m\}\\
v_{i}(t,\cdot)&=u_{i}(t,L-\cdot)\text{ if }i\in\{m+1,...,n\},
\end{split}\right.
\end{equation}
and 
\begin{equation}
\left\{
\begin{split}
\tilde \Lambda_{i}&=\Lambda_{i}\text{ if }i\in\{1,...,m\}\\
\tilde \Lambda_{i}&=-\Lambda_{i}(L-\cdot)\text{ if }i\in\{m+1,...,n\},\\
\tilde B_{i}(\mathbf{v},\cdot)&=B_{i}(\mathbf{u},\cdot)\text{ if }i\in\{1,...,m\}\\
\tilde B_{i}(\mathbf{v},\cdot)&=B_{i}(\mathbf{u},L-\cdot)\text{ if }i\in\{m+1,...,n\}.
\end{split}\right.
\end{equation}
Clearly, $\mathbf{u}$ is a $L^{2}$ solution to the system \eqref{sys1}, \eqref{bound1} if and only if $\mathbf{v}$ is an $L^{2}$ solution to a system of the form \eqref{sys1}, \eqref{bound1} with $\tilde \Lambda$ instead of $\Lambda$ and $\tilde B$ instead of $B$. And now, $\tilde \Lambda$ has only positive components while $\tilde B$ is still Lipschitz with respect to $\mathbf{v}$. Therefore, in this proof, we will assume that $m=n$ and $\Lambda$ has only positive components.
From} \cite[Appendix A.1]{DFP}, $\mathcal{A}+B$ is $\zeta$ dissipative with $\zeta$ independent of $n$ and is a closed operator (in $L^{2}(0,L)$). \textcolor{black}{A definition of an operator $\zeta$ dissipative can be found in \cite[Definition 2.4 and Chapter 5, section 2]{Miyadera}.} Note that \cite{DFP} study systems with a local diagonal source term and positive and constant propagation speeds. However, the proof of these two first points in \cite[Theorem A.1]{DFP} 
only requires $B$ to be Lipschitz
and the propagation speeds 
\textcolor{black}{to be positive and nonvanishing}.
Now we would like to show that $\mathcal{A}+B$ satisfies the following range condition: 
\begin{equation}
\exists \rho_{0}>0,\text{ }\forall\text{ }\rho\in(0,\rho_{0}),\text{ }D(\mathcal{A}+B)\subset Rg(Id-\rho(\mathcal{A}+B)). 
\label{range}
\end{equation}
or equivalently that for any $\mathbf{v}\in D(\mathcal{A}+B),$ there exists $\mathbf{u}\in H^{1}(0,L)$ such that
\begin{equation}
\begin{split}
&\mathbf{u}-\rho\left(\Lambda\partial_{x}\mathbf{u}+B(\mathbf{u},\cdot)\right)=\mathbf{v},\\
&\begin{pmatrix}
\mathbf{u}_{+}(0)\\
\mathbf{u}_{-}(L)
\end{pmatrix}=G\begin{pmatrix}
\mathbf{u}_{+}(L)\\
\mathbf{u}_{-}(0)
\end{pmatrix}.
\label{rangeeq}
\end{split}
\end{equation}
The difficulty comes from the nonlinearity of the equation and this was the main point shown in \cite[Theorem A.1]{DFP}. In our case, all we need to do is to change slightly \textcolor{black}{their} proof to take into account the nonlocal operator and the fact that $\Lambda$ depends on $x$. The latter is easy to take into account by replacing $e^{-\Lambda^{-1} x/\rho}$ by  $e^{-\int_{0}^{x}(\Lambda^{-1}(s)/\rho) ds}$ when integrating, which has a similar behavior \textcolor{black}{(this holds as $\Lambda$ is diagonal)}. To take into account the nonlocal operator, we need to get the estimate \cite[(26)]{DFP} while replacing the 
estimations in \cite[2.2.1]{DFP}, which hold only for the local case when $B$ has a Lipschitz constant independent of $x$. 
But we have from the Lipschitz behavior of $B$ and Cauchy-Schwarz inequality,
\begin{equation}
\begin{split}
&\left|\int_{0}^{x} e^{-\int_{s}^{x}\Lambda^{-1}(v)/\rho dv}\Lambda^{-1}(s)B(\mathbf{u},s) ds\right|\\
&\leq
\lVert B(\mathbf{u},\cdot)\rVert_{L^{2}}\frac{\lVert\Lambda^{-1} \rVert_{\infty}}{2}\rho\int_{0}^{x} \left|e^{\int_{x}^{s}2\Lambda^{-1}(v)/\rho dv}\frac{2\Lambda^{-1}(s)}{\rho}\right| ds\\
&\leq \rho \frac{C_{B}\lVert\Lambda^{-1} \rVert_{\infty}}{2}\left|1-e^{-2\int_{0}^{L}\frac{\Lambda^{-1}(v)}{\rho}} \right| \lVert \mathbf{u}\rVert_{L^{2}}\leq  \rho C_{2} \lVert \mathbf{u}\rVert_{L^{2}}
\end{split}
\end{equation}
where $C_{2}$ is a constant that depends only on the parameters of the system. This enables to recover the estimate \cite[(26)]{DFP} which is then used to apply Arzela-Ascoli Theorem and get the range condition \eqref{range}.
As in Theorem \cite[Theorem A.1]{DFP}, from these three properties ($\zeta$ dissipative, closed operator and range condition) and using \cite[Corollary 5.13 and Remark 2 p.148]{Miyadera}, $\mathcal{A}+B$ generates a \textcolor{black}{nonlinear} semigroup $S$ of type $\zeta$ on $L^{2}(0,L)$ and the Cauchy problem has a unique integral solution $\mathbf{u}\in C^{0}([0,T],L^{2}(0,L))$ (see \cite{Miyadera} for a proper definition of an integral solution). 
\textcolor{black}{Besides, let $\mathbf{u}_{0,n}\in D(\mathcal{A})$, \textcolor{black}{then} from  \cite[Remark 2 p.148]{Miyadera}, the unique integral solution $\mathbf{u}_{n}$ of the Cauchy problem 
with initial condition $\mathbf{u}_{0,n}$ belongs to $C^{1}([0,T]; L^{2}(0,L))\textcolor{black}{\cap C^{0}([0,T]; H^{1}(0,L))}$ and satisfies \eqref{bound1} and \eqref{L2sol}. We can choose a sequence $\mathbf{u}_{0,n}\in D(\mathcal{A})$ such that $\mathbf{u}_{0,n}\rightarrow \mathbf{u}_{0} \in L^{2}(0,L)$, as $D(\mathcal{A})$ is dense in $L^{2}$.} Finally, as $S$ is a semigroup \textcolor{black}{of type $\zeta$} we have (see \cite[Remark p.146]{Miyadera})

\begin{equation}
\lVert S(t)\mathbf{u}_{0}-S(t)\mathbf{u}_{0,n}\rVert_{L^{2}}\leq Ce^{\zeta t}\lVert \mathbf{u}_{0}-\mathbf{u}_{0,n}\rVert_{L^{2}},
\end{equation}
which implies the convergence of $\mathbf{u}_{n}$ to $\mathbf{u}$ in $C^{0}([0,T],L^{2}(0,L))$.
\textcolor{black}{ To conclude we only need to show that this is the unique solution in the sense of Definition \ref{defsol}. Let assume that there is another solution $\mathbf{u}^{(1)}$ with initial condition $\mathbf{u}_{0}$. Let $T>0$. By assumption there exists a sequence $\mathbf{u}_{0,n}^{(1)}\in D(\mathcal{A})$ such that $\mathbf{u}_{n}^{(1)}$ satisfies \eqref{bound1} and \eqref{L2sol} with initial condition $\mathbf{u}_{0,n}^{(1)}$ and $\mathbf{u}_{n}^{(1)}\rightarrow \mathbf{u}^{(1)}$ in $C^{0}([0,T],L^{2}(0,L))$. For any $n\in \mathbb{N}$, $\mathbf{u}_{n}^{(1)}\in C^{1}([0,T],L^{2}(0,L))$, \textcolor{black}{therefore
 $\mathbf{u}_{n}^{(1)}$ is also an integral solution of the Cauchy problem with initial condition $\mathbf{u}_{0,n}^{(1)}$ (see \cite[Remark 2 p.148]{Miyadera}).} Thus, from \cite[Remark p.146]{Miyadera},
\begin{equation}
\lVert \mathbf{u}(t,\cdot)-\mathbf{u}_{n}^{(1)}(t,\cdot)\rVert_{L^{2}}\leq Ce^{\zeta t}\lVert \mathbf{u}_{0}-\mathbf{u}_{0,n}^{(1)}\rVert_{L^{2}},
\end{equation}
and therefore $\mathbf{u}_{n}^{(1)}\rightarrow \mathbf{u}$ in $C^{0}([0,T],L^{2})$, which implies that $\mathbf{u}=\mathbf{u}^{(1)}$ in $C^{0}([0,T],L^{2})$. This holds for any $T>0$, and
ends the proof.}
\end{proof}
\bibliographystyle{elsarticle-num}
\bibliography{Biblio_ISS_2}
\end{document}